\documentclass[preprint,12pt]{elsarticle}

\usepackage{graphicx,amsmath,amsfonts}

\def\bfm#1{\boldsymbol{#1}} 
\def\RR{\mathbb{R}}
\def\NN{\mathbb{N}}

\newtheorem{thm}{Theorem}
\newtheorem{lem}[thm]{Lemma}
\newtheorem{dfn}[thm]{Definition}
\newtheorem{conjecture}[thm]{Conjecture}
\newproof{pf}{Proof}
\newproof{pot}{Proof of Theorem \ref{thm2}}
\journal{Journal of Computational and Applied Mathematics}

\begin{document}

\begin{frontmatter}

\title{A general framework for the optimal approximation of circular arcs by parametric polynomial curves}

\author[address1,address2]{Ale\v{s} Vavpeti\v{c}}
\ead{ales.vavpetic@fmf.uni-lj.si}

\author[address1,address2]{Emil \v{Z}agar\corref{corauth}}
\ead{emil.zagar@fmf.uni-lj.si}
\cortext[corauth]{Corresponding author}

\address[address1]{Faculty of Mathematics and Physics, University of Ljubljana,
Jadranska 19, Ljubljana, Slovenia}
\address[address2]{Institute of Mathematics, Physics and Mechanics, Jadranska 19, Ljubljana, Slovenia}

\begin{abstract}
  We propose a general framework for geometric approximation of circular arcs by parametric
  polynomial curves. The approach is based on constrained uniform approximation of an error
  function by scalar polynomials.
  The system of nonlinear equations for the unknown control points of the approximating polynomial
  given in the B\'ezier form is derived and a detailed analysis provided for some low degree cases
  which might be important in practice.
  At least for these cases the solutions can be, in principal, written in a closed form, and provide the best
  known approximants according to the radial distance. A general conjecture on the optimality of the solution
  is stated and several numerical examples conforming theoretical results are given.
\end{abstract}

\begin{keyword}
geometric interpolation  \sep circular arc \sep parametic polynomial \sep B\'ezier curve
\sep optimal approximation
\MSC[2010] 65D05 \sep  65D07 \sep 65D17
\end{keyword}

\end{frontmatter}


\section{Introduction}
Circular arcs are one of the fundamental geometric primitives and together with straight lines they
have been one of the cornerstones of several graphical and control systems. Later on parametric polynomial
representations of geometric objects have been widely used in applications and successfully upgraded to
non-uniform rational basis splines (NURBS) which nowadays provide an intuitive approach towards to the
construction and modelling of curves and surfaces used in computer aided geometric design (CAGD)
and related fields.
However, there is still an interest in parametric polynomial curves, since they provide even more simple
representations of geometric objects and might still be in use in some software standards. On the other hand,
optimal approximation of special classes of functions or parametric objects by polynomials
has always been a theoretical issue
(Chebyshev alternation theorem \cite{Deutsch-Maserick-SIAM-67},
Stone-Weierstrass approximation theorem \cite{Stone-Weierstrass-Theorem-48}, etc.).
Circular arcs form one such class of curves, since it is well known that a circular arc of positive length can not
be exactly represented in a polynomial form.

A common way to construct parametric polynomial approximants of a circular arc is to
interpolate corresponding geometric quantities. This usually include interpolation of boundary points,
corresponding tangent directions, signed curvatures, etc. The result are so called geometric parametric
polynomial approximants ($G^k$ approximants), which can be put together to geometrically smooth spline curves.

When we are dealing with approximations, the fundamental question is a measure of a distance between
a parametric polynomial approximant and a circular arc. One of the standard measures in this case is the
radial distance measuring the distance of the point on the parametric polynomial
to the corresponding point on the circular arc in the radial direction. It can be shown that under some additional
assumptions it coincides with the well known Hausdorff distance
(\cite{Ahn-Kim-arcs-1997,Jaklic-Kozak-best-circle-preprint}). It is more common to use
a simplified version of the radial distance, the difference between the square of the distance of the point on the
parametric polynomial curve to the center of the circular arc and the square of its radius.
The later one is more attractive since it simplifies the analysis of the
existence and uniqueness of the approximant but still preserves the optimality of the approximation order.
However, it is important to emphasize that the optimal solutions according to this two measures do not
coincide in general.

The list of literature dealing with different types of geometric approximants of a circular arc is long
and we shall mention just the most relevant references according to our approach described later.
Parabolic $G^0$ interpolants were considered in \cite{Morken-91-circles}.
This is actually one of only a few cases where the optimality of the solution was proved.
Different types of
$G^1$ and $G^2$ cubic geometric interpolants were given in early papers
\cite{Dokken-Daehlen-Lyche-Morken-90-CAGD} and \cite{Goldapp-91-CAGD-circle-cubic}.
 Several types of quartic and quintic B\'ezier curves were considered in \cite{Ahn-Kim-arcs-1997},
 and deeper analysis of some geometric quintic approximants can be found in
\cite{Fang-98-CAGD-circle-quintic}. Many new cubic and quartic approximants were also proposed in
\cite{Ahn-Kim-quartic-2007,
Xiaoming_Licai-2010-JCADCG,
Liu-Tan-Chen-Zhang-quartics-2012,
Ahn-Kim-JKSIAM-quartics-13,
Kovac-Zagar-quartics-2014}.
However, in none of the above papers the optimality of the solution has been considered.
The paper \cite{HurKim-2011} is beside \cite{Morken-91-circles} the only one where optimality of some approximants was formally shown. The authors managed to prove it for cubic $G^1$
and quartic $G^2$ approximants.

Some authors also considered the approximation of circular arcs by general degree parametric polynomials.
In \cite{Lyche-Morken-94-Metric},
the Taylor type geometric interpolation, i.e., interpolation at just one point was considered for all odd degree
polynomials. For even degree ones the results can be found in
\cite{Jaklic-Kozak-Krajnc-Zagar-circle-like-07}
and in a more general form in \cite{Jaklic-Kozak-Krajnc-Vitrih-Zagar-conics-13}.
The approximation of the whole circle by Lagrange type approximants can be found
in \cite{Jaklic-circle-CAGD-2016} and in \cite{Jaklic-Kozak-best-circle-preprint}.

The aim of this paper is to present a general framework providing optimal geometric approximants for
general degree $n$ of the parametric polynomial and for any order $k$ of geometric smoothness.
The idea relies on the constrained uniform approximation of the error by scalar polynomials.

The paper is organized as follows. In Section \ref{sec:prelim} the problem is explained in detail and
the radial distance and its simplification are precisely introduced. A general conjecture that the proposed
approach provide optimal solutions is stated. Section \ref{sec:constrained}
concerns constrained uniform approximation of an error function by scalar polynomials. A general theory is
briefly revised and some particular cases needed later are carefully analysed. Next section describes optimal
$G^k$ approximation of circular arcs by parametric polynomial curves. In particular, it provides the
system of nonlinear equations which has to be solved. In Section \ref{sec:cases} some particular
cases are studied in detail. For some of them the optimality is reconfirmed, but for all of them it is shown
that they provide the minimal radial distance among all known approximants. In the last section some
concluding remarks and suggestions for possible future research are given.

\section{Preliminaries}\label{sec:prelim}

We shall consider the following problem. Let $\bfm{c}:[-\varphi,\varphi]\to \RR^2$, $0<\varphi\leq \pi/2$
be a standard nonpolynomial parameterization of a circular arc.
Due to simple affine transformations it is enough to consider the unit circular arcs only, centred at the origin and
symmetric with respect to the first coordinate axis. Thus we can assume that
$\bfm{c}(s)=(\cos s,\sin s)^T$. Our goal is to find
as good as possible approximation of $\bfm{c}$ by parametric polynomial curve $\bfm{p}_n:[-1,1]\to\RR^2$
of degree $n\in\NN$.
It is convenient to express
$\bfm{p}_n=(x_n,y_n)^T$, where $x_n$ and $y_n$ are polynomials of degree at most $n$,
in B\'ezier form, i.e.,
\begin{equation}\label{b_Bern_Bez_form}
  \bfm{p}_n(t)=\sum_{j=0}^n B_j^n(t)\bfm{b}_j,
\end{equation}
where $B_j^n$, $j=0,1,\dots,n$, are (reparameterized) Bernstein polynomials over $[-1,1]$,
given as
\begin{equation*}
  B_j^n(t)=\binom{n}{j}\left(\frac{1+t}{2}\right)^{j}\left(\frac{1-t}{2}\right)^{n-j},
\end{equation*}
and $\bfm{b}_j\in\RR^2$, $j=0,1,\dots,n$, are the control points.\\
The quality of the approximation will be measured by radial distance. For each point
on the parametric curve $\bfm{p}_n$  the closest point on the circular arc
$\bfm{c}$ in the radial direction will be considered. In general, it might happen that no such point exists on
$\bfm{c}$, but some further restrictions on $\bfm{p}_n$ will override this problem.
The formal definition of the radial distance $\widetilde{\psi}_n$  is
\begin{equation*}
  \widetilde{\psi}_n:[-1,1]\to[0,\infty),\quad \widetilde{\psi}_n(t):=\left|\sqrt{x_n(t)^2+y_n(t)^2}-1\right|
  =\left|\|\bfm{p}_n(t)\|_2-1\right|,
\end{equation*}
where $\|\cdot\|_2$ is the standard Euclidean norm on $\RR^2$.
Function $\widetilde{\psi}_n$ is an upper bound for the parametric distance $d_P$, studied in detail in
\cite{Lyche-Morken-94-Metric}. For $\bfm{c}$ and $\bfm{p}_n$ it is defined as
\begin{equation*}\label{def:parametric_distance}
  d_P(\bfm{c},\bfm{p}_n)=\inf_{\rho}\max_{t\in[-1,1]}\left\|(\bfm{c}\circ\rho)(t)-\bfm{p}_n(t)\right\|_2,
\end{equation*}
where $\rho:[-1,1]\to[-\varphi,\varphi]$ is a smooth bijection for which $\rho'>0$.
Clearly, $d_P$ is in general an upper bound for well known Hausdorff distance $d_H$.
If the radial distance between $\bfm{c}$ and $\bfm{p}_n$ is well defined, it can be shown that actually
\begin{equation*}
  d_H(\bfm{c},\bfm{b})=d_P(\bfm{c},\bfm{p}_n)=\max_{t\in [-1,1]}\widetilde{\psi}_n(t)
\end{equation*}
(see \cite{Ahn-Kim-arcs-1997} or \cite{Jaklic-Kozak-best-circle-preprint} for details).
Due to computational reasons it is easier to consider a simplified (signed) radial error
\begin{equation}\label{def:psi_n}
  \psi_n:[-1,1]\to[0,\infty),\quad \psi_n(t):=x_n(t)^2+y_n(t)^2-1=\|\bfm{p}_n(t)\|_2^2-1,
\end{equation}
since no irrational functions are involved but the location of zeros and extrema remains the same as for
$\widetilde{\psi}_n$.
The approximation of circular arc $\bfm{c}$ by parametric polynomial $\bfm{p}_n$
now reduces to the study of optimality of $\psi_n$.

In practice, some additional properties of $\bfm{p}_n$ are required, such as interpolation of boundary points,
tangent directions,\dots More precisely, some geometric interpolation conditions are prescribed at the boundary.
These are given in the following definition.
\begin{dfn}\label{def:Gk}
  A circular arc $\bfm{c}$ and a parametric polynomial curve $\bfm{p}_n$ share a geometric contact of order $k\in\NN$
  at the boundary points $\bfm{c}(\pm\varphi)$,
  if there exists a smooth  regular bijective reparameterization $\rho:[-1,1]\to [-\varphi,\varphi]$ with
  $\rho'>0$, such that
  \begin{equation*}
    \frac{d^j\bfm{p}_n}{dt^j}(\pm 1)=\frac{d^j(\bfm{c}\circ \rho)}{d\tau^j}(\pm 1),
    \quad j=0,1,\dots, k.
  \end{equation*}
  We say that $\bfm{p}_n$ is a $G^k$ approximation of $\bfm{c}$ in this case.
\end{dfn}
The following important result characterizes $G^k$ approximants of circular arcs.
\begin{lem}\label{lem:zerosatboundary}
A parametric polynomial $\bfm{p}_n$ is a $G^k$ approximation of the circular arc $\bfm{c}$
if and only if $\psi_n$ has zeros of multiplicity $k+1$ at $t=\pm 1$.
\end{lem}

The proof of this lemma can be found in \citep{Ahn-Kim-arcs-1997}.
It is well known that parametric polynomials can not reproduce circular arcs of positive length. So for a
$G^k$ approximat $\bfm{p}_n$ of the circular arc $\bfm{c}$
it follows from Lemma \ref{lem:zerosatboundary} that
\begin{equation}\label{def:psi_nk}
  \psi_n(t)=\psi_{n,k}(t):=C\,p_{2n,k}(t),
\end{equation}
where $C\in\RR$ is a nonzero constant and
\begin{equation}\label{def:p2n}
  p_{2n,k}(t):=(1-t^2)^{k+1}q_{n,k}(t)
\end{equation}
is a polynomial of degree $2n$ with $q_{n,k}$  being monic of degree $2n-2k-2$.
By \eqref{b_Bern_Bez_form} and \eqref{def:psi_n}, $C$ and $p_{2n,k}$ both
depend on control points $\bfm{b}_j$, $j=0,1,\dots,n$, which further have to fulfil some additional constraints,
ensuring the $G^k$ continuity from Definition \ref{def:Gk}.
In order to find the best approximant according to \eqref{def:psi_n},
the nonlinear optimization problem have to be solved. If we write
$C=C(\bfm{b}_0,\dots,\bfm{b}_n)$ and $p_{2n,k}(t)=p_{2n,k}(t;\bfm{b}_0,\dots,\bfm{b}_n)$,
then we are looking for
\begin{equation}\label{eq:nonlinear_optimization}
  \min_{\bfm{b}_0,\dots,\bfm{b}_n}\max_{t\in[-1,1]}
  \left|C(\bfm{b}_0,\dots,\bfm{b}_n)\,p_{2n,k}(t;\bfm{b}_0,\dots,\bfm{b}_n)\right|.
\end{equation}
This is definitely very hard nonlinear constrained optimization task.
The authors in several papers simplified it in a way
that they have chosen a polynomial $q_{n,k}$ from \eqref{def:p2n} in advance and then
minimized the constant $C$.
This can be done, e.g.,
by prescribing zeros of $q_{n,k}$. However, the quality of the approximant heavily relies on the selection of
zeros and optimality is  not guaranteed. The only known direct optimizations
\eqref{eq:nonlinear_optimization} seem to be in \cite{Morken-91-circles} and in \cite{HurKim-2011},
where the authors considered an optimal quadratic $G^0$, cubic $G^1$ and quartic $G^2$ approximation of
circular arcs. All these problems were dealing with one parametric families of approximants, and it seems that there are no results known about optimal approximants when several parameters are involved.

Here we propose a new general framework which might provide optimal approximants in any case.
We again choose $q_{n,k}$ of degree $2n-2k-2$, but now in a way that it provides a minimum of
\begin{equation}\label{eq:best_monic_min}
  \|p_{2n,k}\|=
  \max_{t\in[-1,1]}\left|(1-t^2)^{k+1}\,q_{n,k}(t)\right|.
\end{equation}
The polynomial $p_{2n,k}$, which minimizes \eqref{eq:best_monic_min}, will be denoted by
$p_{2n,k}^*$, and the corresponding $q_{n,k}$ by $q_{n,k}^*$.
The idea comes from the constrained uniform approximation of zero function on $[-1,1]$
by monic polynomials and will be considered in detail in the next section. However, this choice does not
a priori guarantee the optimality of the approximant as one might quickly conclude from uniform polynomial
approximation of functions. There might exist approximants which do not provide minimal
$\|p_{2n,k}\|$,
but they provide a constant $C$ small enough that corresponding $\left|\psi_{n,k}\right|$ would be smaller
than the one arising from $p_{2n,k}^*$.
However, results of the present paper show that there is some hope this actually can not happen.

It is clear that $p_{2n,k}^*$, which minimizes \eqref{eq:best_monic_min}, does not depend
on $\bfm{b}_j$, $j=0,1,\dots,n$. It depends only on $n$, $k$, and the properties of the norm
\eqref{eq:best_monic_min}.
Once $p_{2n,k}^*$ is determined,
then control points $\bfm{b}_j$, $j=0,1,\dots,n$, are given as a solution of the system of nonlinear
equations, and we are left with the minimization of $C(\bfm{b}_0,\dots,\bfm{b}_n)$ (i.e., we have to choose
a solution providing minimal $\left|C(\bfm{b}_0,\dots,\bfm{b}_n)\right|$). The main purpose of this paper
is to show that the proposed approach reproduces the above mentioned optimal solutions obtained
in \cite{Morken-91-circles} and in \cite{HurKim-2011}
and provides new solutions for $G^0$ cubic and $G^1$ quartic approximants possessing the smallest known error.
This leads us to the following conjecture.
\begin{conjecture}\label{conj:optimal}
  The best $G^k$ geometric approximant $\bfm{p}_n$ of the circular arc $\bfm{c}$
  according to the error measure $\psi_{n,k}$, given by \eqref{def:psi_nk},
  arises from the choice $p_{2n,k}^*$
  determining $\bfm{b}_j$, $j=0,1,\dots,n$, with minimal $\left|C(\bfm{b}_0,\dots,\bfm{b}_n)\right|$.
\end{conjecture}
In the following section the general approach to the construction of $p_{2n,k}^*$ will be described.

\section{Constrained uniform approximation}\label{sec:constrained}
In this section the optimal approximation of the zero function by polynomials of
the form $p_{2n,k}(t)=(1-t^2)^{k+1}\,q_{n,k}(t)$, where $q_{n,k}$ is a monic polynomial of degree
$2n-2k-2$, will be considered.
In particular, we shall study the following problem: For any $k,n\in\NN$, such that
$0\leq k<n$, find a monic polynomial $q_{n,k}^*$ of degree $2n-2k-2$ for which
$p_{2n,k}^*(t)=(1-t^2)^{k+1}q_{n,k}^*(t)$
has minimal max norm on $[-1,1]$.\\
Suppose that $q_{n,k}\colon [-1,1]\to\RR$ is a monic polynomial of degree $2n-2k-2$.
Let us define the polynomial
$p_{2n,k}$ of degree $2n$ by
\begin{equation}\label{def:p}
  p_{2n,k}(t)=(1-t^2)^{k+1}q_{n,k}(t).
\end{equation}
It follows from  \citep{Loeb-Moursund-Schumaker-Taylor-69} that there exists
the unique monic polynomial $q_{n,k}^*$ of degree $2n-2k-2$, such that
$p_{2n,k}^*(t)=(1-t^2)^{k+1}q_{n,k}^*(t)$ has minimal max norm
over all polynomials of the form \eqref{def:p}.
The polynomial $p_{2n,k}^*$ is characterized by the following property
\citep[Theorem 3.1]{Loeb-Moursund-Schumaker-Taylor-69}:
There exist $2n-2k-1$ points $-1< a_0<\ldots<a_{2n-2k-2}<1$ such that
$\|p_{2n,k}^*\|=\left|p_{2n,k}^*(a_0)\right|$ and
$p_{2n,k}^*(a_0)=(-1)^i p_{2n,k}^*(a_i)$ for $i=1,\dots,2n-2k-2$.
Since we are dealing with polynomials defined over symmetric interval $[-1,1]$,
and since for every monic polynomial $r$ of even degree also $t\mapsto\tfrac 1 2(r(t)+r(-t))$ is a monic
of norm no greater than $\|r\|$, the polynomials $q_{n,k}^*$ and $p_{2n,k}^*$ must be even.
Hence
\begin{equation}\label{eq:p2nform}
  p_{2n,k}^*(t)=(1-t^2)^{k+1}(t^2-t_1^2)\cdots (t^2-t_{n-k-1}^2)
\end{equation}
 for some $0<t_1<\ldots<t_{n-k-1}<1$.
Some special cases which can be analyzed analytically and will be needed later for the construction
of particular geometric approximants,  will now be considered in detail.

\subsection{The case $k=0$}
The polynomial $p^*_{2n,0}$ has exactly $2n$ single roots. Two of them are on the boundary of the interval.
From a general theory of uniform approximation by polynomials it follows that
$p^*_{2n,0}$ is  a scaled and dilated Chebyshev polynomial $T_{2n}$ of degree $2n$, more precisely
\begin{equation}\label{eq:p2n0}
  p_{2n,0}^*(t)=-\frac{2^{1-2n}}{\cos^{2n}\left(\tfrac{\pi}{4n}\right)}
  T_{2n}\left(\cos\left(\frac{\pi}{4n}\right)t\right).
\end{equation}
It is easy to deduce $q_{n,0}^*$ since its zeros must be precisely the $2n-2$ interior zeros of
$p_{2n,0}^*$. In particular, if $n=2$, we have $t_1=\sqrt{3-2\sqrt{2}}$, and for $n=3$ the pozitive zeros are
$t_1=2-\sqrt{3}$ and $t_2=\sqrt{3}-1$.

\subsection{The case $k=n-2$}
In this case the derivative of $p_{2n,n-2}^*$ can be written as
$$
  \frac{dp^{*}_{2n,n-2}}{dt}(t)=2n(1-t^2)^{n-2}t(1-a-t^2)
$$
for some $a\in(0,1)$. Direct integration then leads to
\begin{equation*}
  p_{n,n-2}^*(t)=(1-t^2)^{n-1}(t^2-1+\tfrac{n}{n-1}a).
\end{equation*}
The polynomial $p_{2n,n-2}^*$ has its extrema at
$-\sqrt{1-a}$, $0$, and $\sqrt{1-a}$. By the characterization of the best approximant
the equality $p_{2n,n-2}^*(\sqrt{1-a})=-p_{2n,n-2}^*(0)$ has to be fulfilled . Hence, in order
to find $p_{2n,n-2}^*$, the equation
\begin{equation}\label{eq:f1}
  \varphi(a):=a^n+n\,a-(n-1)=0
\end{equation}
has to be solved and the solution must be in $I_n:=(0,1-1/n)$. Since $\varphi(0)=1-n<0$ and
$\varphi(1-1/n)=(1-1/n)^n>0$, the solution is indeed in $I_n$ and it is unique due to the Descartes rule of
signs (\cite{Albert-43}).
In particular, for $n=3$ we have $a=\sqrt[3]{\sqrt{2}+1}-\sqrt[3]{\sqrt{2}-1}$ and
\begin{equation}\label{eq:bestp61}
  p_{6,1}^*(t)=(1-t^2)^2\left(t^2-1+\frac{3}{2}a\right).
\end{equation}
Thus $t_1=\sqrt{1-\frac{3}{2}a}$.

\subsection{The case $k=n-3$}
Similarly as in the previous subsection, we can write
\begin{equation*}
  \frac{dp_{2n,n-3}^*}{dt}(t)=2n(1-t^2)^{n-3}t(1-a-t^2)(1-b-t^2),\quad 0<b<a<1.
\end{equation*}
Integration of the previous form gives
$$
  p_{2n,n-3}^*(t)=(1-t^2)^{n-2}\left(-(1-t^2)^2+\frac{n}{n-1}(a+b)(1-t^2)
  -\frac{n}{n-2}ab\right).
$$
Since the extrema of $p_{2n,n-3}^*$ are $\pm\sqrt{1-a}$, $\pm\sqrt{1-b}$, and $0$, characterization of
the best approximat implies
\begin{equation}\label{eq:condp2}
  p_{2n,n-3}^*(0)=-p_{2n,n-3}^*(\sqrt{1-a})=p_{2n,n-3}^*(\sqrt{1-b}).
\end{equation}
The last equality in \eqref{eq:condp2} leads to
\begin{align*}
-a^{n-1}\left(\tfrac{n}{(n-1)(n-2)}b-\tfrac1{n-1}a\right)=b^{n-1}\left(\tfrac{n}{(n-1)(n-2)}a-\tfrac1{n-1}b\right).
\end{align*}
Multiplication by $\tfrac {n-1} {b^n}$ gives
\begin{equation}\label{eq:fi2}
  \zeta(\lambda):=\lambda^n-\frac{n}{n-2}\lambda^{n-1}-\frac{n}{n-2}\lambda+1=0,
\end{equation}
where $\lambda:=\tfrac b a\in(0,1)$.
Since $\zeta(0)=1>0$ and $\zeta(1)=-4/(n-2)$, $\zeta$ must have at least one root on $(0,1)$.
But then by the Descartes rule of signs there must be exactly two positive roots. Due to the symmetry, roots of
$\zeta$ must appear in pairs $\lambda$,$1/\lambda$. Consequently we have the unique root on $(0,1)$.
In particular, if $n=4$, we have $\zeta(\lambda)=\lambda^4-2\lambda^3-2\lambda+1=0$ and
\begin{equation}\label{lambda4}
  \lambda=(\sqrt{3}-\sqrt{2}\sqrt[4]{3}+1)/2.
\end{equation} Since $b=\lambda a$, the first equality in \eqref{eq:condp2}
and the fact that $\zeta(\lambda)=0$ imply an equation for unknown $a$, namely
\begin{equation}\label{eq:a4}
\left(a-\frac{1}{\lambda}\right)^2
\left(\lambda^2(2\lambda-1)\,a^2+2\lambda(2\lambda-1)\,a+6\lambda^3+6\lambda-3\right)=0.
\end{equation}
However,  $\lambda<1$, and $a$ must be the unique positive zero of the second factor in \eqref{eq:a4}.
Some calculations reveal that $a$ can be written as
\begin{equation}\label{eq:asol}
  a=\sqrt{1+\sqrt{3}+\sqrt{24+14\sqrt{3}}}-\frac{1}{2}\left(1+\sqrt{3}+\sqrt{2}\sqrt[4]{3}\right).
\end{equation}
The polynomial $p_{8,1}^*$ then reads as
\begin{equation}\label{eq:p81}
  p_{8,1}^*(t)=-(1-t^2)^2(t^4+\frac{2}{3}(2(\lambda+1)a-3)t^2+\frac{1}{3}(3-4(\lambda+1)a+6\lambda a^2)),
\end{equation}
where $\lambda$ is given by \eqref{lambda4} and $a$ by \eqref{eq:asol}.
Two positive zeros $t_{1,2}$ of $p_{8,1}^*$ can be found as a solution of the quadratic equation arising
from the quartic factor in \eqref{eq:p81}.

\section{Optimal $G^k$ approximation of circular arcs}\label{sec:Gkapprox}

In previous section constrained uniform minimization by polynomials was studied in detail.
In order to use the obtained results, let us consider a general problem of $G^k$ approximation of circular arcs by
parametric polynomials of arbitrary degree $n$. 	Suppose that the approximant $\bfm{p}_n$ is given by
\eqref{b_Bern_Bez_form}. Quite clearly, $\bfm{p}_n$ has $2n+2$ free parameters, i.e., the coordinates
of the control points $\bfm{b}_j$, $j=0,1,\dots,n$. Since the circular arc $\bfm{c}$ is symmetric
with respect to abscissa, so is the approximant $\bfm{p}_n$. Consequently, its control points must be symmetric
too, and the number of free parameters reduces to $n+1$. Additionally,  $G^0$ condition
at a particular point prescribes two parameters (one control point), and each $G^\ell$ condition,
$1\leq \ell\leq k$, reduces the number of free parameters by one (\cite{Farin-97-CAGD}).
Finally, $G^k$ approximant
is determined by $n+1-2-k=n-k-1$ parameters. Particularly,
if $k=0$, the first and the last control points
must be $\bfm{b}_0=(\cos{\varphi},-\sin{\varphi})^T$ and $\bfm{b}_n=(\cos{\varphi},\sin{\varphi})^T$.
If $k=1$, additionally $\bfm{b}_1=\bfm{b}_0+d\,\bfm{c}'(-\varphi)$ and
$\bfm{b}_{n-1}=\bfm{b}_n-d\,\bfm{c}'(\varphi)$ for some $d>0$. Some similar, but more complicated
relations can be derived for $k\geq 2$, too.

Since we are interested in optimal $G^k$ approximation, we shall follow our proposed approach and choose
\begin{equation}\label{error_psi_condition}
\psi_{n,k}(t)=\|\bfm{p}_n(t)\|_2^2-1=C\,p_{2n,k}^*(t)=C\,(1-t^2)^{k+1}\,q_{n,k}^*(t),
\end{equation}
where $p_{2n,k}^*$ minimizes \eqref{eq:best_monic_min}. By \eqref{eq:p2nform},
$q_{n,k}^*$ has $2n-2k-2$ symmetric roots on $(-1,1)$. Let $0<t_1<t_2<\cdots<t_{n-k-1}<1$
be the positive ones. Then
\begin{equation}\label{p_system}
  \psi_{n,k}(t_i)=\|\bfm{p}_n(t_i)\|_2^2-1=0,\quad i=1,2,\dots,n-k-1,
\end{equation}
is a system of $n-k-1$ nonlinear equations for the $n-1-k$ unknown parameters
determining the approximant $\bfm{p}_n$. It might have several solutions,
and we are interested in that one which implies the minimal absolute value of $C$ in
\eqref{error_psi_condition}. Since due to the symmetry $\psi_{n,k}$ must have an extreme point at $0$, and
all the extrema are by construction of the same magnitude, we have
\begin{equation}\label{eq_C}
  C=\frac{\|\bfm{p}_n(0)\|_2^2-1}{q_{n,k}^*(0)}.
\end{equation}
Among all possible solutions of \eqref{p_system}, we thus choose the one providing $\bfm{p}_n$ for which
\eqref{eq_C} has the minimal absolute value.

There is a little hope that the problem can be solved for general $k$ and $n$. Thus we will concentrate on some
specific low degree cases and we shall confirm that our approach either reproduces the best solution or implies
the solution for which the error is the smallest among all by now known approximants. This gives some hope that the
Conjecture \ref{conj:optimal} might hold true for any $k$ and $n$.

\section{Some particular cases}\label{sec:cases}

The first nontrivial case is $n=2$ and $k=0$, i.e., quadratic approximation via interpolation of
boundary points of the circular arc $\bfm{c}$. This case has already been considered in \cite{Morken-91-circles}
and the best solution was characterized. We reconsider it here to demonstrate the elegance of our approach
and to prove that our conjecture holds true in this case. We also derive an asymptotic expansion of $C$, which was not provided in \cite{Morken-91-circles}.

\begin{figure}[htb]
\centering
\includegraphics[width=0.9\textwidth]{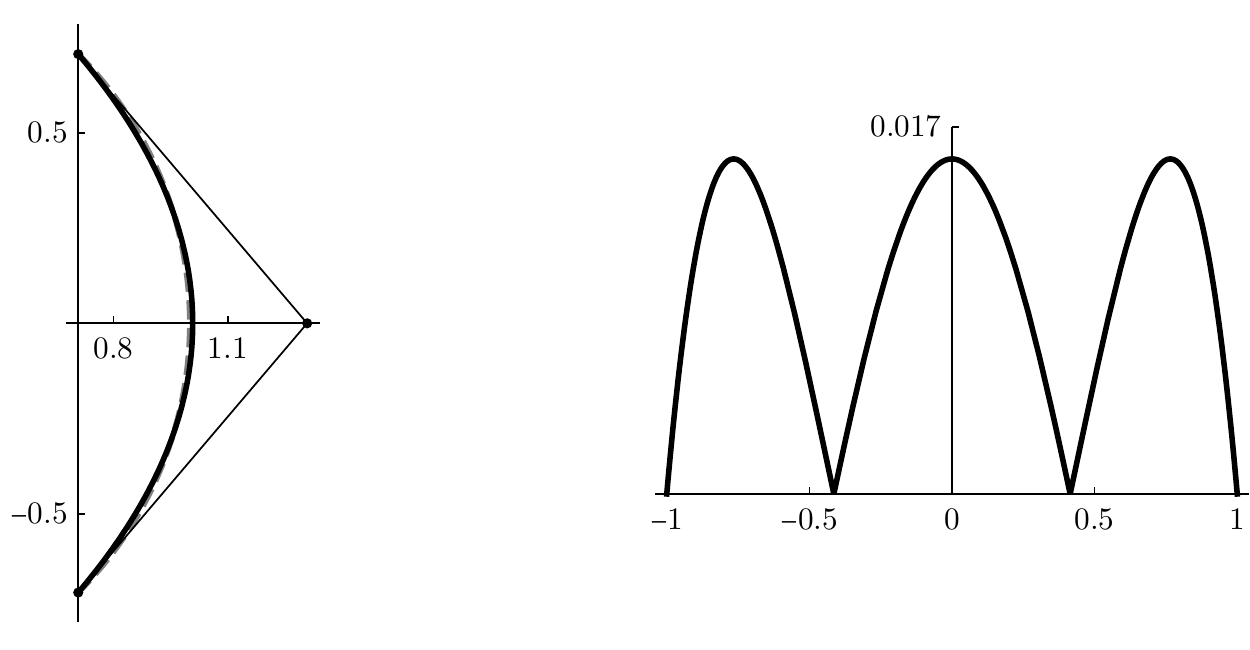}
\caption{Left: Circular arc given by central angle of $\pi/2$ (dashed gray) together with the optimal
  quadratic $G^0$ approximant (black). Right: The error $\left|\psi_{2,0}^*\right|$.}
  \label{fig:quadraticG0}
\end{figure}

The approximant $\bfm{p}_2$ is determined by three control points
\begin{equation*}
  \bfm{b}_0=(\cos\varphi,-\sin\varphi)^T,\quad
  \bfm{b}_1=(\xi,0)^T \quad {\rm and}  \quad
  \bfm{b}_2=(\cos\varphi,\sin\varphi)^T,
\end{equation*}
where $\xi>0$. By \eqref{eq:p2n0} the only positive zero of $p_{4,0}^*$ is
$t_1=\sqrt{3-2\sqrt{2}}$ and by \eqref{p_system} the unknown parameter $\xi$
must fulfil the equation $\|\bfm{p}_2(t_1)\|_2^2-1=0$, or equivalently
\begin{equation*}
  (6-4\sqrt{2})\,\xi^2 +
   (12\sqrt{2}-16 )\cos{\varphi}\,\xi + (3-2\sqrt{2})\cos{2\varphi}+7-6\sqrt{2}=0.
\end{equation*}
Since $\xi$ must be positive, the only admissible solution is
\begin{equation*}
  \xi=-\sqrt{2}\cos\varphi + \sqrt{2 + 2\sqrt{2} + \cos^2\varphi}
\end{equation*}
and coincides with the solution derived in \cite{Morken-91-circles}.
From \eqref{eq_C} it is quite easy to get an asymptotic expansions $C=-\varphi^4/4+{\cal O}(\varphi^5)$ and
\begin{equation*}
  \max_{t\in[-1,1]}\left|\psi_{2,0}(t)\right|=\frac{3-2\sqrt{2}}{4}\varphi^4+{\cal O}(\varphi^5)\approx
  0.0429\varphi^4+{\cal O}(\varphi^5).
\end{equation*}
This also proves that the  $G^0$ quadratic
approximant is fourth order accurate. Some easy calculations further reveal the Hausdorff distance
$d_H(\bfm{c},\bfm{b})\approx 0.0214\,\varphi^4+{\cal O}(\varphi^5)$.
An example of the best circular $G^0$ approximant and its error are shown on Fig.~\ref{fig:quadraticG0}.

The case $n=2$ and $k=1$ is not interesting, since the inner control point $\bfm{b}_1$ is uniquely
determined by $G^1$ condition as the intersection of tangent lines to the boundary points of the circular arc.
Thus, there is nothing to be optimized.

\begin{figure}[hbt]
  \centering
    \includegraphics[width=0.6 \textwidth]{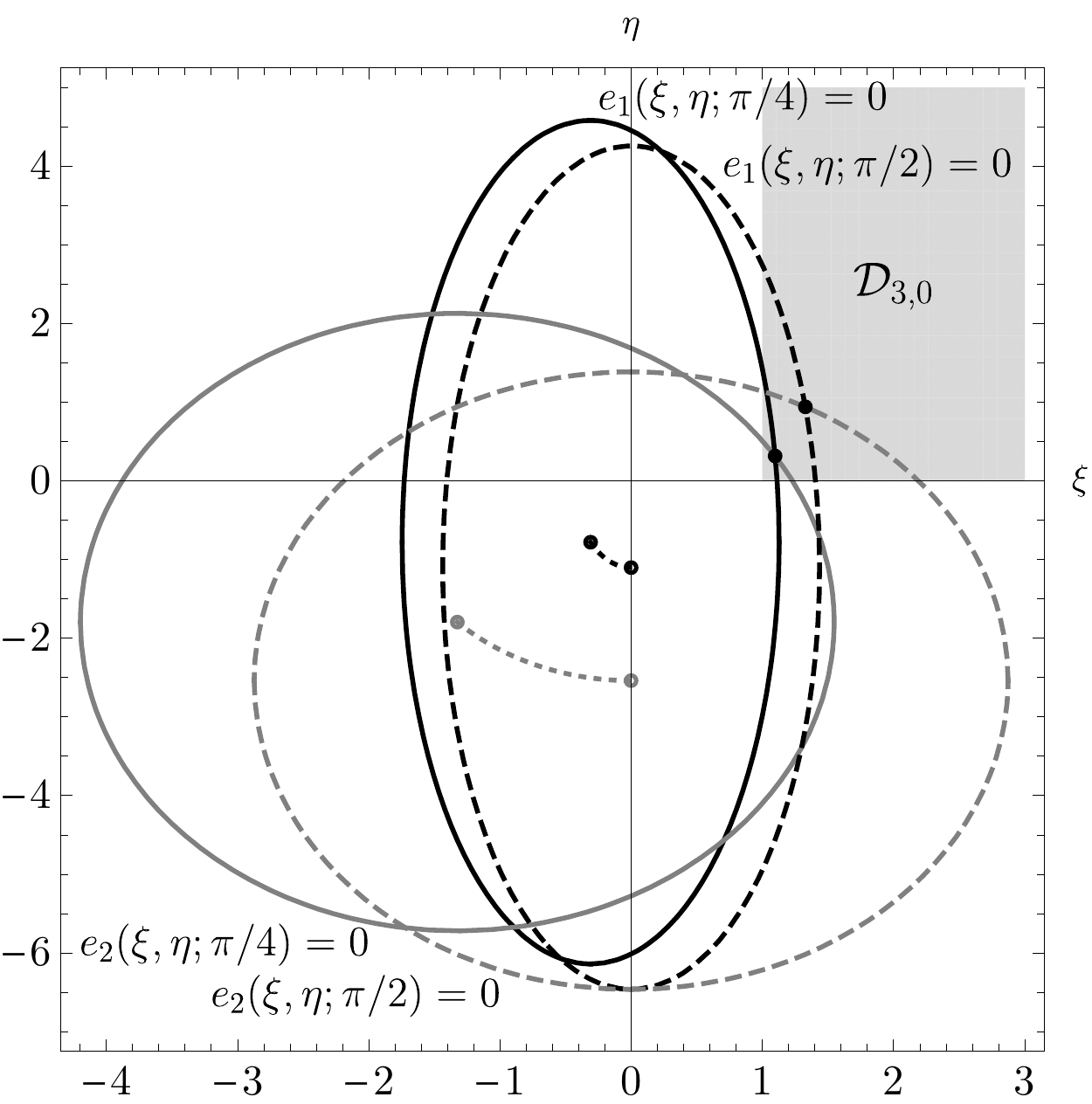}
  \caption{Graphs of ellipses $e_1(\xi,\eta;\pi/4)=0$ (solid black), $e_1(\xi,\eta;\pi/2)=0$ (dashed black),
  $e_2(\xi,\eta;\pi/4)=0$ (solid gray) and $e_2(\xi,\eta;\pi/2)=0$ (dashed gray). In the middle of the figure
  are loci along which centres of ellipses are moving for $\varphi\in[\pi/4,\pi/2]$.
  The region ${\mathcal D}_{3,0}$
  is shadowed upper right. Solutions of the nonlinear system \eqref{eq:ellipses} for $\varphi=\pi/4$ and
  $\varphi=\pi/2$ are black points in  ${\mathcal D}_{3,0}$.}
  \label{fig:ellipses}
\end{figure}

Let us consider the case $n=3$ and $k=0$ now. It was partially considered in
\cite{Goldapp-91-CAGD-circle-cubic}. The author reported that it leads to the solution of the nonlinear biquadratic
system, but no proof of the existence and uniqueness of the solution was provided. Here we fill this gap by the
formal proof arising from our approach. Due to $G^0$ conditions, the control points of the approximant
$\bfm{p}_3$ must be
\begin{equation*}
  \bfm{b}_0=(\cos\varphi,-\sin\varphi)^T,\quad \bfm{b}_1=(\xi,-\eta)^T,\quad
  \bfm{b}_2=(\xi,\eta)^T,\quad \bfm{b}_3=(\cos\varphi,\sin\varphi)^T,
\end{equation*}
where obviously $\xi>1$ and $\eta>0$.
By \eqref{eq:p2n0}
\begin{equation*}
  p_{6,0}(t)=-\frac{1}{32\cos\left(\frac{\pi}{12}\right)^6}
  T_6\left(\cos\left(\frac{\pi}{12}\right)t\right),
\end{equation*}
and its positive roots on $(0,1)$ are $t_1=2-\sqrt{3}$ and $t_2=\sqrt{3}-1$.
The system of nonlinear equations for $\xi$ and $\eta$
\begin{equation}\label{system30}
  f_i(\xi,\eta;\varphi):=\psi_{3,0}(t_i)=0,\quad i=1,2,
\end{equation}
now follows from \eqref{p_system}.
Some further computations reveal that the equations \eqref{system30} actually represent two ellipses. More
precisely, \eqref{system30} is equivalent to
\begin{equation}\label{eq:ellipses}
  e_i(\xi,\eta;\varphi)=\frac{\left(\xi-p_i(\varphi)\right)^2}{a_i^2}
  +\frac{\left(\eta-q_i(\varphi)\right)^2}{b_i^2}-1=0,\quad i=1,2,
\end{equation}
where the coordinates of the centres of the ellipses are
\begin{align}
  (p_1(\varphi),q_1(\varphi))&=
  \left(\frac{1}{9} \left(3-4 \sqrt{3}\right)\cos\varphi,
  \frac{1}{9} \left(-3-4 \sqrt{3}\right)\sin\varphi\right),
  \label{center1}\\
  (p_2(\varphi),q_2(\varphi))&=
  \left(\frac{1}{9} \left(-3-8 \sqrt{3}\right)\cos\varphi,
  \frac{1}{9} \left(-9-8 \sqrt{3}\right)\sin\varphi\right),
  \label{center2}
\end{align}
and the semiaxes are
\begin{align*}
  (a_1,b_1)&=\left(\frac{2}{9}\left(3+2\sqrt{3}\right),\frac{2}{9}\left(12+7\sqrt{3}\right)\right),
  \\
  (a_2,b_2)&=\left(\frac{4}{9}\left(3+2\sqrt{3}\right),\frac{2}{9}\left(9+5\sqrt{3}\right)\right).
  \\
\end{align*}
Thus the solution of the system of nonlinear equations \eqref{system30} is given by the
intersection of two ellipses \eqref{eq:ellipses}. Therefore, it is enough to show that
this two particular ellipses intersect in ${\cal D}_{3,0}:=\{(\xi,\eta);\ \xi>1,\ \eta>0\}$.
An example of such ellipses for $\varphi=\pi/4$ and $\varphi=\pi/2$ is shown on Fig. \ref{fig:ellipses}.
It is clearly seen that there is precisely one intersection in ${\cal D}_{3,0}$. This will now be confirmed
by the following lemma.

\begin{figure}[hbt]
  \centering
    \includegraphics[width=0.9 \textwidth]{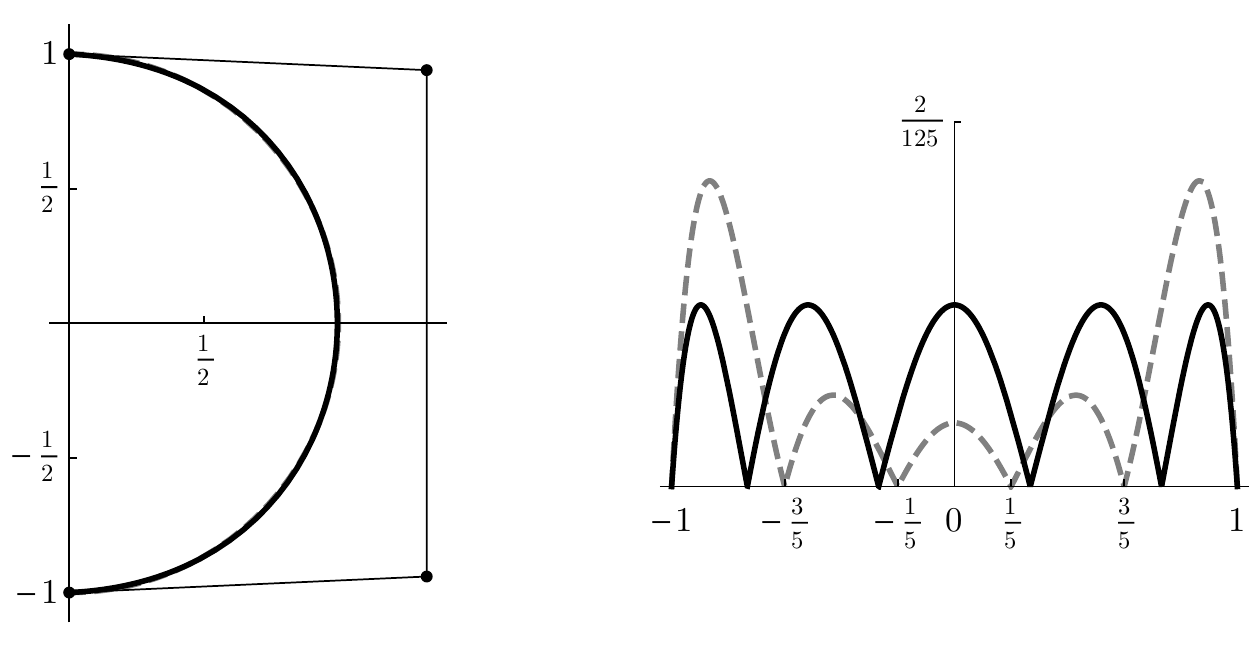}
  \caption{Left: Circular arc given by central angle of $\pi$ (dashed gray, invisible) together with the optimal
  cubic $G^0$ approximant (black). Right: Graphs of errors $\left|\psi_{3,0}^*\right|$ and
  $\left|\psi_{3,0}\right|$ possessing  equidistant zeros.}
  \label{fig:cubicG0}
\end{figure}

\begin{lem}\label{lem:intersection_ellipses}
  The ellipses $e_i(\xi,\eta;\varphi)=0$, $i=1,2$, intersect precisely at one point
  in ${\cal D}_{3,0}$.
\end{lem}
\begin{pf}
  The prove will base on the fact which is depicted in Fig.~\ref{fig:ellipses} in the shadowed region.
  We shall prove that similar situation appears for each $\varphi\in(0,\pi/2)$.
  A closer look to ellipses $e_i=0$, $i=1,2$, reveals that they actually rotate when $\varphi$ runs
  over an angle interval $(0,\pi/2)$. It can also be shown that their centres rotate along two
  other ellipses which are uniquely determined by \eqref{center1} and \eqref{center2}.
  Let $\xi_i$ be the solutions of
  $e_i(\xi,0;\varphi)=0$ on the boundary of ${\mathcal D}_{3,0}$, $i=1,2$, respectively.
  Similarly, let $\eta_i$ be the solutions of $e_i(1,\eta;\varphi)=0$  on the boundary of
  ${\mathcal D}_{3,0}$,
  $i=1,2$, respectively. To prove that there is precisely one intersection
  point in ${\cal D}_{3,0}$, it is enough to see that
  \begin{equation}\label{relation_xi_eta}
   (\xi_1-\xi_2)(\eta_1-\eta_2)<0.
  \end{equation}
  Since $e_i(\xi,0;\varphi)=0$,  $i=1,2$,
  are quadratic equations in $\xi$, the intersections $\xi_i$, $i=1,2$, are easily determined.
  It turns out that
  \begin{align*}
    \xi_1&=\frac{1}{9}\left(3-4\sqrt{3}\right)\cos\varphi+\sqrt{\frac{1}{54}\left(19+52\sqrt{3}\right)
    +\frac{1}{54}\left(37-20\sqrt{3}\right)\cos 2\varphi},\\
    \xi_2&=-\frac{1}{9}\left(3+8\sqrt{3}\right)\cos\varphi+\sqrt{\frac{1}{27}\left(74+59\sqrt{3}\right)
    +\frac{1}{27}\left(38+5\sqrt{3}\right)\cos 2\varphi}.
  \end{align*}
  It is similarly easy to see that the solutions of $e_i(1,\eta;\varphi)=0$, $i=1,2$, are
  \begin{align*}
    \eta_1&=-\frac{1}{9}\left(3+4\sqrt{3}\right)\sin\varphi\\
    &+\sqrt{\frac{2}{27}\left(199+116\sqrt{3}\right)+
    \frac{2}{27}\left(37+20\sqrt{3}\right)\cos\varphi}\,\sin\frac{\varphi}{2},\\
    \eta_2&=-\frac{1}{9}\left(9+8\sqrt{3}\right)\sin\varphi\\
    &+\sqrt{\frac{1}{27}\left(362+213\sqrt{3}\right)+
    \frac{1}{27}\left(182+99\sqrt{3}\right)\cos\varphi}\,\sin\frac{\varphi}{2}.
  \end{align*}
  A straightforward computation using some basic properties of trigonometric functions
  leads us to \eqref{relation_xi_eta}.\qed
\end{pf}
Thus we can state the following theorem.
\begin{thm}\label{thm:G30}
  The system of nonlinear equations \eqref{system30} has a unique admissible solution $(\xi^*,\eta^*)$
  in ${\cal D}_{3,0}$. Consequently, there exists the unique cubic $G^0$ approximant of the circular arc,
  given by the central angle $\varphi\in(0,\pi/2]$.
  The error of the approximation is
  $$
    \max_{t\in[-1,1]}\left|\psi_{3,0}^*(t)\right|=\frac{1}{64}\left(26-15\sqrt{3}\right)\varphi^6
    +{\cal O}(\varphi^8)\approx 0.0003 \varphi^6+{\cal O}(\varphi^8).
  $$
\end{thm}
\begin{pf}
    The existence and the uniqueness of the optimal solution follows from previous lemma.
    The asymptotic expansion can be obtained from the Taylor expansion of $\psi_{3,0}^*(0)$
    around $\varphi=0$ considering analytic solution for $\xi^*$ and
    $\eta^*$ (rather longish formulae which will not be written here) and using some computer
    algebra system.\qed
\end{pf}
Although the exact formulae for $\xi^*$ and $\eta^*$ from Theorem \ref{thm:G30}
can, in principal, be obtained, they will probably be evaluated numerically in practice.
One can use a particular iterative method (e.g., Newthon-Raphson method),
since quite accurate starting points for the iteration can be obtained by finding an
approximate intersection of the ellipses in ${\mathcal D}_{3,0}$.
For $\varphi=\pi/2$, the optimal
solution becomes particularly simple, namely $\xi^*=4\sqrt{2+4\sqrt{3}}/9$ and
$\eta^*=(5+2\sqrt{3})/9$. Fig.~\ref{fig:cubicG0} shows the approximant together with the error
in this case.
As a comparison, we took the approximant arising from the choice of $p_{6,0}$ having uniformly distributed
zeros on $[-1,1]$, i.e., $\pm 1$, $\pm 1/5$ and $\pm 3/5$. The corresponding error is shown on
Fig.~\ref{fig:cubicG0} and is much bigger.

For the next case let us consider $G^1$ approximation by cubics. As it was already mentioned, the problem
was considered in \cite{HurKim-2011}. The authors characterized an optimal solution and they have
show its optimality.
One of the reasons that they manage to prove the optimality is the fact that the family of approximants depend
on only one parameter, the case already observed in $G^0$ quadratic case. We shall see that our approach
again simplifies the analysis significantly.

\begin{figure}[hbt]
  \centering
    \includegraphics[width=0.9 \textwidth]{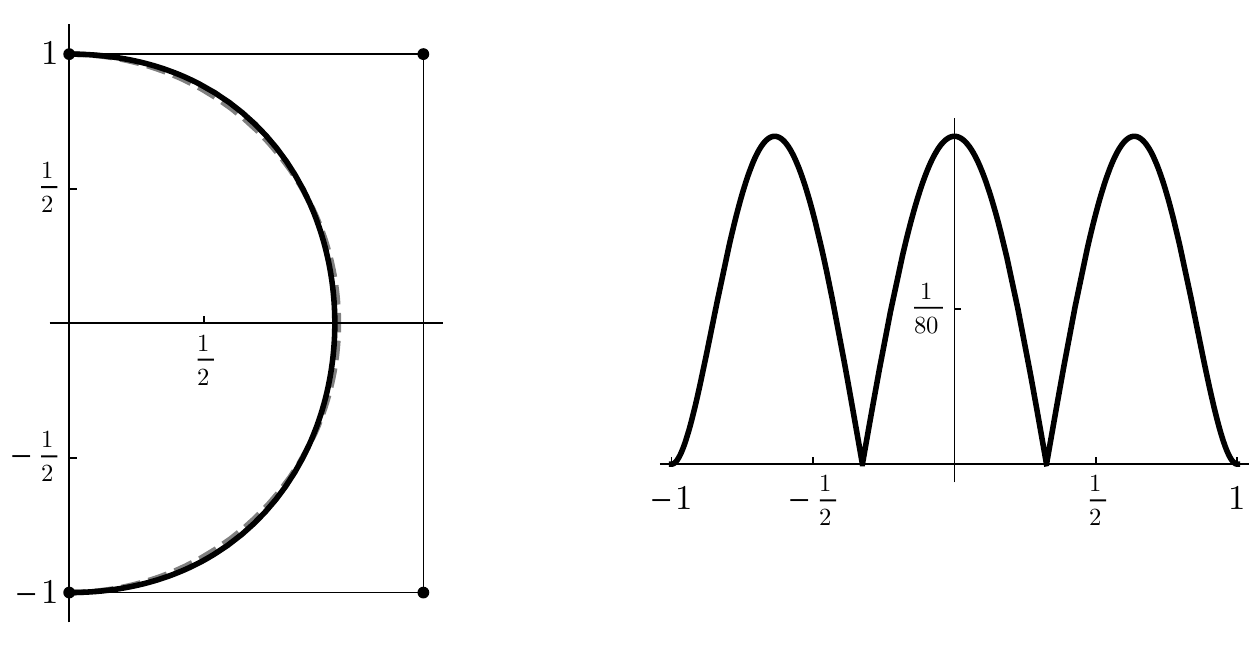}
  \caption{Left: Circular arc given by central angle of $\pi$ (dashed gray) together with the optimal
  cubic $G^1$ approximant (black). Right: Graph of the error  $\left|\psi_{3,1}^*\right|$.}
  \label{fig:G31}
\end{figure}

Control points of the $G^1$ cubic approximant $\bfm{p}_3$ can be now given as one parameter family of the form
\begin{align*}
  &\bfm{b}_0=(\cos\varphi,-\sin\varphi)^T, &&\bfm{b}_1=\bfm{b}_0+d\,(\sin\varphi,\cos\varphi)^T,\\
  &\bfm{b}_3=(\cos\varphi,\sin\varphi)^T,&&\bfm{b}_2=\bfm{b}_3-d\,(-\sin\varphi,\cos\varphi)^T,
\end{align*}
where $d>0$. By \eqref{eq:bestp61}, the only positive zero of $p_{6,1}^*$ is $t_1=\sqrt{1-\frac{3}{2}a}$,
where $a=\sqrt[3]{\sqrt{2}+1}-\sqrt[3]{\sqrt{2}-1}$.
According to \eqref{p_system}, we only have to solve one quadratic
equation for $d$, namely
\begin{align}
f(d;\varphi)&:=(9(b^2-1)(1+\cos 2\varphi)+12b)\,d^2\nonumber\\
&-4\sin 2\varphi(3b^2-2b-3)\,d+8\sin^2\varphi((b-1)^2-2)=0,\label{eq:d}
\end{align}
where $b=\sqrt[3]{\sqrt{2}-1}$. It is easy to see that $f(0;\varphi)<0$ and $f''(d;\varphi)>0$
for $\varphi\in(0,\pi/2)$, which proves that \eqref{eq:d} has a unique solution on $(0,\infty)$
for any $\varphi\in(0,\pi/2)$.
It is also easy to see that this solution coincides with the optimal solution obtained in \cite{HurKim-2011}.
The solution for $\varphi=\pi/2$ again significantly simplifies to $d=\sqrt{\frac{2}{3}\left(\frac{1}{b}-b+2\right)}$.
The approximant for this case together with the error is shown in Fig. \ref{fig:G31}.

Finally, let us consider the quartic $G^1$ case for which the optimality of the solution has not been
studied yet. The control points of the parametric polynomial approximant
$\bfm{p}_4$ are
\begin{align*}
  &\bfm{b}_0=(\cos\varphi,-\sin\varphi)^T,\quad \bfm{b}_1=\bfm{b}_0+d\,(\sin\varphi,\cos\varphi)^T,
  \quad \bfm{b}_2=(\xi,0),\\
  &\bfm{b}_4=(\cos\varphi,\sin\varphi)^T,\quad \bfm{b}_3=\bfm{b}_4-d\,(-\sin\varphi,\cos\varphi)^T,
\end{align*}
where $\xi$ and $d$ should be in
${\mathcal D}_{4,1}:=\left\{(\xi,d);\ \xi>1,\ d>0\right\}$.  Quartic $G^1$ approximant again form two parameter family as it was the case
in $G^0$ cubic approximation. It is thus expected that the problem is much harder as the $G^1$ cubic
case, and this is probably also the reason that there is no result on the optimal solution in the literature.
We will again follow our approach and show that the solution provides the smallest known error.

Let $t_{1,2}\in(0,1)$ be positive zeros of the polynomial $p_{8,1}^*$, defined by \eqref{eq:p81}.
Note that this two zeros can be given in a closed form, since we only have to solve a quadratic equation with
exactly known coefficients determined by \eqref{eq:p81}. The nonlinear system for the unknowns $\xi$ and $d$
is given by \eqref{p_system} as
\begin{equation}\label{sys:g12}
  g_i(\xi,d;\varphi):=\psi_{4,1}(t_i)=0,\quad i=1,2.
\end{equation}
It can be shown again that  $g_i(\xi,d;\varphi)=0$ represent (possibly) degenerated ellipses, but this time
their configuration is more complicated as in $G^0$ cubic case. However, analytic representation of the ellipses
(however quite complicated and given by longish formulae) allow us to use the same conclusions as in the $G^0$ cubic
case. Let $d_i$ be the solutions of the (quadratic) equations $g_i(1,d;\varphi)=0$,
and $\xi_i$ be the solutions of the (quadratic) equations $g_i(\xi,0;\varphi)=0$, $i=1,2$,
on the boundary of ${\mathcal D}_{4,1}$. The existence of the unique solution
of the nonlinear system \eqref{sys:g12} is confirmed if the relation, similar to \eqref{relation_xi_eta}, namely
\begin{equation}\label{rel:dxi}
  (d_1-d_2)(\xi_1-\xi_2)<0,
\end{equation}
is fulfilled. Unfortunately, numerical examples show that the above relation might not be true for some
small $\varphi$ (definitely smaller that $\varphi=\pi/4$). But in this case there are always two solutions,
and one of them provides much smaller radial error as the other one. The analysis of this fact is beyond of the
scope of this paper, and we will show the existence of the unique solution for $\varphi=\pi/4$ only. Furthermore,
we will show that the error of the approximant is the smallest one among the errors of quartic $G^1$ approximants
studied in the literature.

\begin{figure}[htb]
\centering
\begin{minipage}{0.48\textwidth}
\includegraphics[width=1\textwidth]{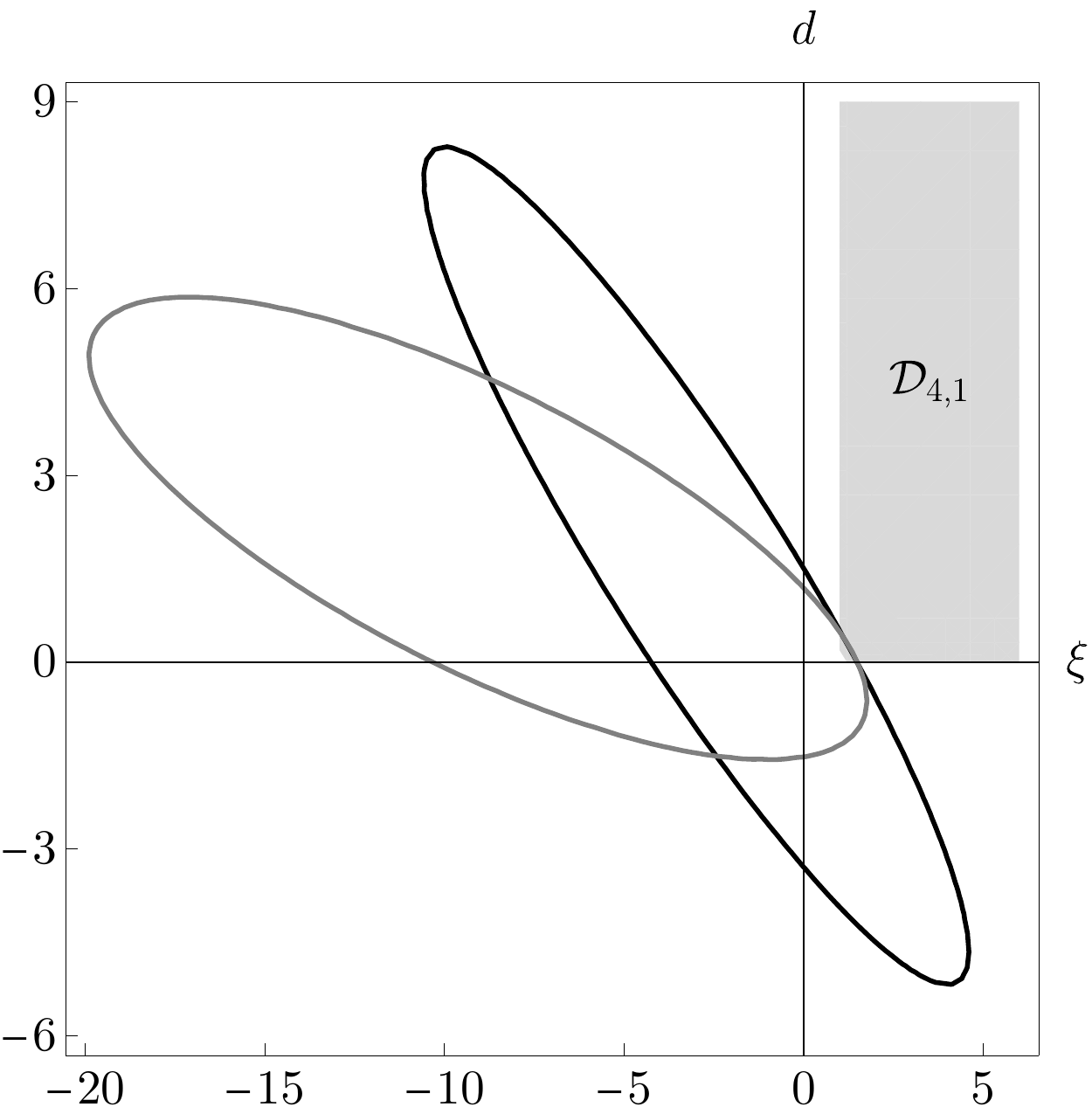}
\end{minipage}
\begin{minipage}{0.48\textwidth}
\includegraphics[width=1\textwidth]{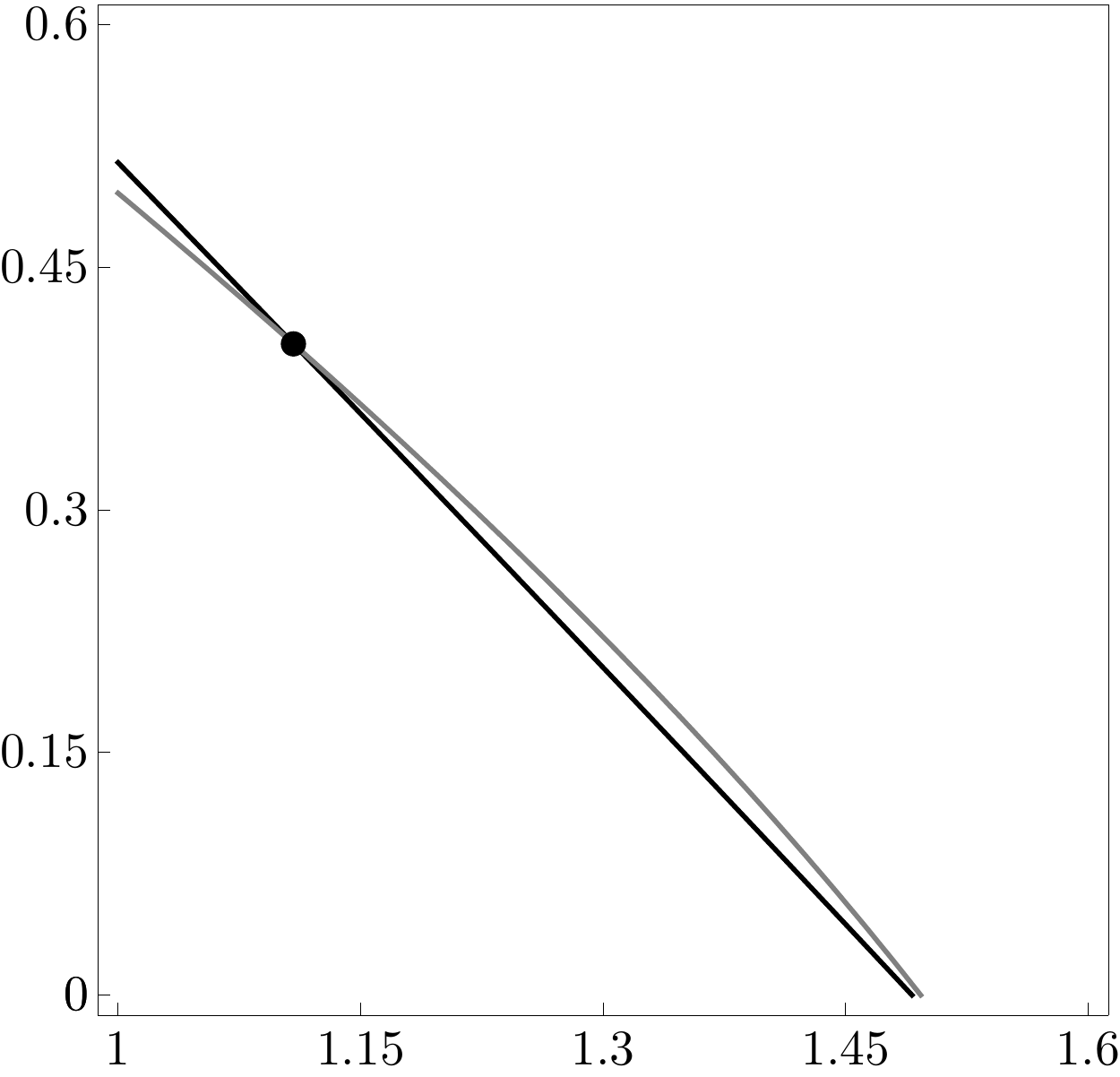}
\end{minipage}\\[4mm]
\caption{Left: Graphs of the ellipses $g_i(\xi,d;\pi/4)=0$, $i=1,2$. The upper right corner is the admissible region
${\mathcal D}_{4,1}$. Right: The zoom of the arcs of the ellipses in the region ${\mathcal D}_{4,1}$.}
\label{fig:G41ellipses}
\end{figure}

If $\varphi=\pi/4$, then the solutions $\xi_i$ and $d_i$, $i=1,2$, can be found analytically by solving
quadratic equations with exact coefficients. It is then easy to find numerical values $d_1=0.514871$,
$d_2=0.495957$, $\xi_1=1.49096$ and $\xi_2=1.496410$ and \eqref{rel:dxi} follows. This proves the existence
of the unique admissible quartic $G^1$ approximant of the circular arc given by the inner angle $\pi/2$.
The ellipses defined by $g_1$ and $g_2$ are shown on
Fig. \ref{fig:G41ellipses}. The same figure also shows the behaviour of the elliptic arcs in the region
 ${\mathcal D}_{4,1}$.
In Table \ref{table:G41} Hausdorff distances of the known quartic $G^1$ approximants and the circular arc given by the
inner angle $\pi/2$ are collected. Clearly our solution provides smaller error as the existing ones.

\begin{table}[htb!]
\begin{center}
\begin{tabular}{|c|r|c|}
\hline
Zeros of $p_{81}$  & Hausdorff distance  & Reference\\ \hline
$-1,-1,-1,-1,1,1,1,1$   & $3.50\times10^{-5}$    & \cite{Ahn-Kim-arcs-1997}\\ \hline
$-1,-1,-w_1,-w_1,w_1,w_1,1,1$ & $4.75\times 10^{-6}$ & \cite{Kovac-Zagar-quartics-2014}
\\ \hline
$-1,-1,-1,0,0,1,1,1$ & $3.55\times 10^{-6}$ & \cite{Ahn-Kim-arcs-1997}\\ \hline
$-1,-1,0,0,0,0,1,1$ & $2.03\times 10^{-6}$ & \cite{Ahn-Kim-CAD-quartics-2007}\\ \hline
$-1,-1,-1/2,0,0,1/2,1,1$ & $1.11\times 10^{-6}$ & \cite{Liu-Tan-Chen-Zhang-quartics-2012}\\ \hline
$-1,-1,-3/5,0,0,3/5,1,1$ & $1.08\times 10^{-6}$ & \cite{Xiaoming_Licai-2010-JCADCG}\\ \hline
$-1,-1,-z_1,0,0,z_1,1,1$ & $7.60\times 10^{-7}$ & \cite{Ahn-Kim-JKSIAM-quartics-13}\\ \hline
$-1,-1,-t_2,-t_1,t_1,t_2,1,1$ & $6.34\times 10^{-7}$ & This paper.\\ \hline
\end{tabular}
\caption{Hausdorff distances of known quartic $G^1$ approximants to the circular arc given by inner angle
of $\pi/2$. Zeros in the table are: $w_1=\sqrt{2}-1$,
$z_1=\frac{1}{3}\sqrt{6-4\sqrt{3}+2\sqrt{6}\sqrt{\sqrt{3}-1}}$ and $0<t_1<t_2<1$ are
zeros of $p_{8,1}^*$.}
\label{table:G41}
\end{center}
\end{table}

\begin{figure}[htb]
\centering
\includegraphics[width=0.85\textwidth]{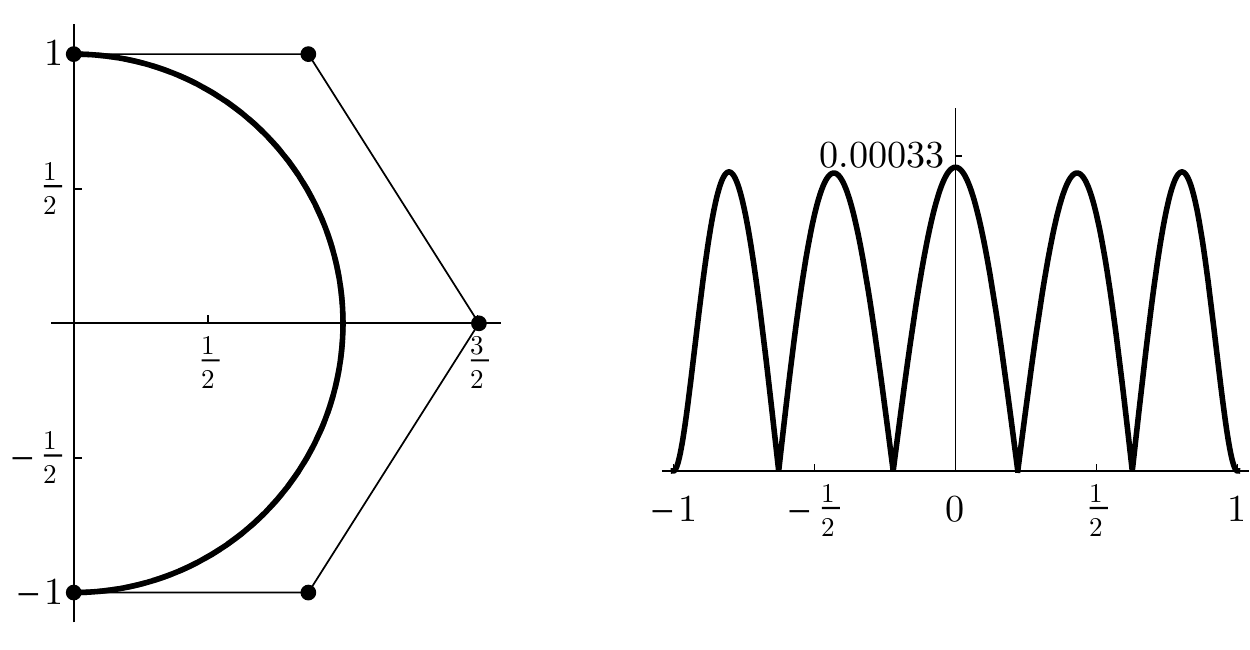}
\caption{Left: Circular arc given by central angle of $\pi$ (dashed gray, invisible) together with the optimal
  quartic $G^1$ approximant (black). Right: Graph of the error  $\left|\psi_{4,1}^*\right|$.}
\label{fig:G41semicircle}
\end{figure}

The same analysis could be done for any fixed angle $\pi/4<\varphi\leq\pi/2$. However, the proof for a general
angle seems to be much more complicated and relies on the powerful computer algebra system. As an example,
the solution for $\varphi=\pi/2$ is shown on Fig. \ref{fig:G41semicircle}. Numerical values for the
parameters are $\xi=1.50506$ and $d=0.87152$.
These values can be obtained also in a closed form, but the expressions are to long to be presented here.

It is also possible to study $G^2$ quartic approximation using our approach. The family of geometric
approximants depends again on just one parameter and the analysis simplifies significantly comparing to
\cite{HurKim-2011}. Only one scalar nonlinear equation has to be analysed with the same solution
as it was provided in \cite{HurKim-2011}.
\section{Conclusion}

We have presented a new approach to the solution of optimal geometric approximation of circular arc by
parametric polynomial curves. It is based on constrained uniform approximation by polynomials.
The solutions obtained by the proposed approach coincide with known optimal approximants
in $G^0$ quadratic, $G^1$ cubic and $G^2$ quartic case. A general conjecture on the optimality of the
geometric approximants was stated and some particular cases which have been studied in detail are confirming it.
As a future work the proposed approach can be used for some other low degree geometric approximants,
but some higher order algebraic equations as \eqref{eq:f1} or \eqref{eq:fi2} have to be solved first.
It is not to be expected that the solutions can be given in radicals and numerical procedures are unavoidable.
But even more important issue would be the proof of Conjecture \ref{conj:optimal}. This would assure that
our approach gives the best possible approximants according to the simplified radial distance measure.
However, this does not solve the problem of optimality in the sense of Hausdorff distance. To do this,
one has to consider radial error $\widetilde{\psi}$ as a measure of an error. Numerical experiments
indicate that this is much more difficult problem.
\vspace{1cm}

{\noindent \sl Acknowledgments.}{The first author acknowledges the financial support from the Slovenian Research Agency (ARRS) by research core funding no. P1-0292, J1-7025, J1-8131.
The second author was supported in part by the program P1-0288 and the grant J1-7256 from ARRS.}


\end{document}